\documentclass{article}

\usepackage{authblk}
\usepackage{amsmath,amsfonts}

\def\Theta{\Phi}
\def\cz{\mathbb{C}}
\def\im{\text{\tt i}}

\begin{document}

\title{Multivariate exponential analysis from the minimal number of samples}

\author{Annie Cuyt and Wen-shin Lee\\
Department of Mathematics and Computer Science\\ University of Antwerp (CMI)\\ 
Middelheimlaan 1, B-2020 Antwerpen, Belgium
\\ \texttt{\{annie.cuyt,wen-shin.lee\}@uantwerpen.be}}
\maketitle  

\begin{abstract}
The problem of multivariate exponential analysis or sparse interpolation
has received a lot of attention, especially with respect to the number of
samples required to solve it unambiguously. In this paper we 
show how to bring the
number of samples down to the absolute minimum of $(d+1)n$ where $d$ is
the dimension of the problem
and $n$ is the number of exponential terms. 
To this end
we present a fundamentally different approach for the multivariate
problem statement. We combine a one-dimensional exponential analysis method such as
ESPRIT, MUSIC, the matrix pencil or any Prony-like method, with some linear systems of equations
because the multivariate exponents are
inner products and thus linear expressions in the parameters. 
\\


\noindent\textbf{Mathematics Subject Classication (2010):} 42B99, 42A15. 
\end{abstract}

\section{Introduction}

Multivariate exponential analysis is a classical problem at the basis of
many application domains (such as, for instance,
\cite{650355,Man2000,YILMAZER2006796,4244728}) 
that recently has gained a lot of attention.
The problem statement is that of recovering the vectors $\phi_j \in \cz^d,
j=1, \ldots, n$ and the coefficients $\alpha_j \in \cz, j=1, \ldots, n$ in the $d$-variate 
$n$-sparse sum
\begin{multline*}
f(x) := f(x_1, \ldots, x_d) = \sum_{j=1}^n \alpha_j \exp\left( \langle
\phi_j, x \rangle \right), \\ 
x= (x_1, \ldots, x_d), \qquad \phi_j=(\phi_{j1}, \ldots, \phi_{jd}), \qquad
\langle \phi_j, x \rangle= \sum_{i=1}^d \phi_{ji} x_i,
\end{multline*}
from $(d+1)n$ samples of $f(x_1, \ldots, x_d)$, which is the minimal number
of samples because it equals the number of parameters in the problem
statement. 

When $d=1$ then the problem can be solved using a variety of Prony-based
algorithms \cite{Be.Ti:det:88,Sc:mul:86,Ro.Ka:esp:89,Hu.Sa:mat:90}, 
in which the identification of the $\phi_j$ and $\alpha_j$ is
separated and taken care of in two stages. 
The frequencies $\phi_j, j=1, \ldots, n$ are obtained from a 
generalized eigenvalue or polynomial rooting problem, while the linear 
coefficients $\alpha_j, j=1, \ldots, n$ are computed from a Vandermonde
system of linear equations \cite[pp.~378--382]{Hi:int:56}.
Input to these algorithms are $2n$ samples of
$f(x)$ at some equidistant points $f(s\Delta), s=0, \ldots, 2n-1$. This
number of samples is minimal if $n$ is known. Otherwise at least one
more sample is required to identify the sparsity $n$. For more details on
the latter we refer to \cite{Ka.Le:ear:03,Cu.Le:spa:16}.

Several computational methods were developed to solve the problem also when
$d>1$, from straightforward generalizations to more sophisticated
approaches, all of them using more than a minimum of $(d+1)n$ samples though.
It should be obvious to the reader that the challenge is not to
recover inner products $\langle \phi_j, x \rangle$ and the associated
coefficients $\alpha_j$ for $j=1, \ldots, n$,
from $2n$ equidistant samples in 
higher-dimensional space. Under modest conditions this 
can be achieved using the univariate techniques mentioned above. 
Instead, the challenge is to recover
the individual $\phi_{ji}, j=1, \ldots, n, i=1, \ldots, d$ and the
coefficients $\alpha_j$. We describe the state of the art in multivariate
exponential analysis
and explain how our approach differs from it.

The one-dimensional 
matrix pencil method was generalized to the 2-dimen\-sio\-nal
matrix enhancement and matrix pencil method (MEMP) \cite{Hu:est:92} 
and can be extended to higher dimensions in a straightforward manner. 
It uses a Hankel-block-Hankel matrix to decompose the 2-dimensional
problem into two one-dimensional problems reflecting each dimension. This
decomposition introduces an additional challenge though, namely that of matching
or pairing the information computed from the one-dimensional problems
\cite{Ro.Na:est:01}. Moreover,
when constructing a uniform $d$-dimensional grid of sample points, the
amount of information is $O(n^d)$. 

Solving the problem along some one-dimensional subspace, 
in other words computing some projection 
such as in \cite{Pl.Wi:how:13,Po.Ta:par:13*b} 
requires only $O(n)$ samples. 
Using an adaptive sampling scheme and
under some mild condition on the coefficients, this remains valid in the
2-dimensional case \cite{Wischerhoff2016}.
However, in \cite{Di.Is:par:15} is shown that there is no finite set of 
(independently of $f$) predefined lines
for which the bivariate reconstruction problem has a unique
solution. A lower bound for the number of samples in the reconstruction 
when $d=2$ is $O(n^2)$. In order to solve the pairing problem,
\cite{Di.Is:par:15} reformulates the problem as a non-convex optimization
problem, which is not computationally feasible for practical purposes.

Rather than projecting on one-dimensional subspaces, a symbolic approach
based on \cite{Sauer2017}
is developed in \cite{Sa:pro:16} making use of constructive ideal theory
and multivariate polynomial interpolation. The largest number of required 
samples in
this setting is estimated to be $O((d+1)n^2 \log^{2d-2}n)$. In the same
corner one finds \cite{Ku.Pe.ea:mul:16} and \cite{Pe.Pl.ea:rec:} which
obtain the multivariate exponents as common roots of a finite system of
$d$-variate polynomials. Still making use of $O(n^d)$ samples however,
algebraic geometry theory now guarantees the correct pairing and recovery.

The method we propose differs significantly from all of the above, not
only in its informational usage which can be as low as $(d+1)n$, but also 
in its approach which only makes use of a 1-dimensional Prony technique
combined with some linear systems of equations
because the individual $\phi_{ji}$
appear linearly in the $\langle \phi_j, x\rangle$.
The presented multivariate exponential analysis technique results from
ideas that were initially formulated in
\cite{Cu.Le:sma:12,Cu.Le:sma:12*b}: a so-called identification shift in
the sampling strategy allows to overcome any ambiguity in the exponential
analysis.

After this state of the art of the literature, 
Sections 2 and 3 deal with the ideal case where some mild assumptions
are verified and only $(d+1)n$
evaluations are necessary, thus generalizing Prony's result where $2n$
samples solve a univariate exponential analysis problem. 
In Section 4 the most general case is detailed, requiring slightly more
samples because the assumptions do not hold. 
An analysis of the worst case scenario and an algorithm for the detection
of $n$ is presented in Section 5. Finally,
the new algorithm is illustrated with an example in Section 6.

\section{Multivariate exponential analysis}

As surveyed in the introduction, up to now computational methods 
require more samples than the minimal number, for one or other reason.
%
%
We now explain how the problem statement can also be solved in the
multivariate setting using the minimal number $(d+1)n$ of samples. The trick
to achieve this is to split the set of samples in two subsets, namely $2n$
equidistant samples and another $(d-1)n$ samples that may but need not be
equidistant in the higher-dimensional space
(they cannot be entirely unstructured though). 
We discuss the use of the $2n$ equidistant samples in this
section and that of the additional $(d-1)n$ samples in Section 3.
For now we assume in the multivariate setting that the value of $n$ is known.
How to detect $n$ is further discussed in Section 5. 

Let $\Delta=(\Delta_1, \ldots, \Delta_d) \not=(0, \ldots, 0)$ and 
$|\Im \phi_{ji}| < \pi/|\Delta_i|, j=1, \ldots, n, i=1, \ldots, d$
\cite{Ny:cer:28,Sh:com:49}, where the function $\im$ returns the imaginary part
of a complex number.
Let us 
sample $f(x_1, \ldots, x_d)$ at the points $s\Delta, s=0, \ldots, 2n-1$:
\begin{equation}
F_{s} := f(s\Delta_1, \ldots, s\Delta_d), \qquad s=0, \ldots, 2n-1.
\label{1stbatch}
\end{equation}
For the time being, we also assume that the sampling direction $\Delta$ is such
that the values $\exp(\langle \phi_j, \Delta \rangle)$, $j=1, \ldots, n$
are mutually distinct. How to deal with collisions in these values is
described in Section 4.

Following the univariate scheme \cite[pp.~378--382]{Hi:int:56} 
the coefficients $\beta_i, i=0, \ldots$, $n-1$ of the polynomial
\begin{equation}
B(z) = \prod_{j=1}^n \left( z-\exp \left( \langle \phi_j, \Delta \rangle 
\right) \right) = z^n + \beta_{n-1} z^{n-1} + \ldots + \beta_0 \label{pol}
\end{equation}
can be obtained from the $n\times n$ Hankel system of linear equations
\begin{equation}
\begin{pmatrix} F_{0} & F_{1} & \cdots & F_{n-1} \\ F_{1} & \cdots
& & F_{n} \\
\vdots & & & \vdots \\ F_{n-1} & F_{n} & \cdots & F_{2n-2} 
\end{pmatrix}
\begin{pmatrix} \beta_0 \\ \vdots \\ \beta_{n-1} \end{pmatrix} = 
- \begin{pmatrix} F_{n} \\ \vdots \\ F_{2n-1} \end{pmatrix}, \label{HA}
\end{equation}
or the roots $\exp(\langle \phi_j, \Delta \rangle), j=1, \ldots,n$
of $B(z)$ can be found as the generalized eigenvalues $\lambda$ of the problem 
\begin{equation}
\begin{pmatrix} F_{1} & F_{2} & \cdots & F_{n} \\ F_{2} & \cdots
& & F_{n+1} \\
\vdots & & & \vdots \\ F_{n} & F_{n+1} & \cdots & F_{2n-1}
\end{pmatrix} v = \lambda
\begin{pmatrix} F_{0} & F_{1} & \cdots & F_{n-1} \\ F_{1} & \cdots
& & F_{n} \\
\vdots & & & \vdots \\ F_{n-1} & F_{n} & \cdots & F_{2n-2}
\end{pmatrix} v, \qquad v \in \cz^n. \label{EV}
\end{equation}
So we can recover the expressions $\exp(\Phi_j)$ where
\begin{equation}
\Phi_j = \langle \phi_j, \Delta \rangle, \qquad j=1, \ldots, n.
\label{Phi}
\end{equation}
Although we have not yet identified the individual $\phi_{ji}, j=1, \ldots,
n, i=1, \ldots, d$,
nothing prevents us from already computing the linear coefficients
$\alpha_j$ from one of the $n \times n$ Vandermonde systems
\begin{multline}
\begin{pmatrix} 
\exp(k\Phi_1) & \exp(k\Phi_2) & \cdots
& \exp(k\Phi_n) \\ \vdots & & & \vdots \\ \exp((k+n-1)\Phi_1) &
\exp((k+n-1)\Phi_2) & \cdots & \exp((k+n-1)\Phi_n) \end{pmatrix} 
\begin{pmatrix} \alpha_1 \\ \vdots \\ \alpha_n \end{pmatrix} \\ = 
\begin{pmatrix} F_{k} \\ \vdots \\ F_{k+n-1} \end{pmatrix},
\qquad 0 \le k \le n. \label{VM}
\end{multline}
The latter can also be replaced by the $2n \times n$ Vandermonde system 
involving all samples, which is then
solved in the least squares sense, as recommended in the case of real-life and hence 
noisy data.

\section{Identification shifts}

In order to extract the $\phi_{ji}, j=1, \ldots, n, i=1, \ldots, d$ from the
$\Phi_j, j=1, \ldots, n$, still under the assumption that the values
$\exp(\Phi_j), j=1, \ldots, n$ are mutually distinct,
some additional samples are required. For this purpose we choose a set
$\{\Delta, \delta_1, \ldots, \delta_{d-1}\}$ of $d$ linearly independent
vectors in $\cz^d$. The additional samples are taken along
a linear combination of $\Delta$ and some $\delta_i, i=1, \ldots, d-1$: 
\begin{multline}
F_{\ell i} := f(\kappa_{\ell i}\Delta + \delta_i)
=f(\kappa_{\ell i}\Delta_1 + \delta_{i1}, \ldots, 
\kappa_{\ell i}\Delta_d + \delta_{id}), \\ \ell=1, \ldots, n, 
\quad i=1, \ldots, d-1 \label{IDsamples}
\end{multline}
where the $\kappa_{\ell i}, \ell=1, \ldots, n$ 
for fixed $i$ are taken to be mutually distinct.
A simple choice for $\kappa_{\ell i}$ for all $i$
is $\kappa_{\ell i}= \ell-1$. Then the
additional samples are taken equidistantly
along independent shifts $\delta_i$ with
respect to the original vector $\Delta$, in other words $F_{\ell i}=
f((\ell-1)\Delta+\delta_i)$. 
At the same time we assume that
$$|\Im \langle \phi_j, \delta_i/||\delta_i|| \rangle | <
\pi/||\delta_i||, \qquad j=1, \ldots, n, \quad i=1, \ldots, d$$
in order to comply with the Shannon-Nyquist conditions formulated in 
\cite{Ny:cer:28,Sh:com:49}.
We call these vectors $\delta_i, i=1, \ldots, d-1$ identification shifts for
reasons that will become apparent: they allow to identify the individual
$\phi_{ji}$ from the computed $\Phi_j$. For this identification we exploit
the fact that the $\phi_{ji}$ appear linearly in the $\Phi_j$ and hence we
turn our attention to systems of linear equations rather than to 
multivariate polynomial root solving or structured 
generalized eigenvalue problems. 

Consider for fixed $i=1, \ldots, d-1$, meaning for a chosen linearly
independent shift vector 
$\delta_i$, the following Vandermonde-like system of linear equations:
\begin{equation}
\begin{pmatrix} \exp(\kappa_{1i} \Phi_1) & \exp(\kappa_{1i}\Phi_2) & \cdots & 
\exp(\kappa_{1i}\Phi_n) \\ \\ \exp(\kappa_{2i}\Phi_1) & \cdots & &
\exp(\kappa_{2i}\Phi_n) \\ \vdots & & & \vdots \\ \\ \exp(\kappa_{ni}\Phi_1)
& \cdots & & \exp(\kappa_{ni}\Phi_n)) \end{pmatrix}
\begin{pmatrix} A_{1i} \\ \vdots \\ A_{ni} \end{pmatrix} = \begin{pmatrix}
F_{1i} \\ \vdots \\ F_{ni} \end{pmatrix}. \label{VMlike}
\end{equation}
Since we know $\exp(\Phi_j), j=1, \ldots, n$ and have chosen $\kappa_{\ell
i}, \ell=1, \ldots, n$, with $i$ fixed,
the Vandermonde-like coefficient matrix can easily be composed. 
Note that for the choice $\kappa_{\ell i}=\ell-1$ the Vandermonde-like
coefficient matrix coincides with that of \eqref{VM} where $k=0$. The unknowns 
$A_{ji}, j=1, \ldots, n$ come from a reinterpretation of the samples
$F_{\ell i}$ as
$$F_{\ell i} = f(\kappa_{\ell i}\Delta+\delta_i)=
\sum_{j=1}^n \alpha_j \exp\left( \langle \phi_j, \delta_i \rangle) \right)
\exp\left( \langle \phi_j, \kappa_{\ell i}\Delta \rangle \right),
\qquad \ell=1, \ldots, n$$
with 
$$A_{ji} = \alpha_j \exp\left( \langle \phi_j, \delta_i \rangle) \right),
\qquad j=1, \ldots, n$$
and 
$$\exp(\kappa_{\ell i}\Phi_j) =
\exp\left( \langle \phi_j, \kappa_{\ell i}\Delta \rangle \right), 
\qquad \ell, j=1, \ldots, n.$$
The values $A_{ji}/\alpha_j$ equal
$${A_{ji} \over \alpha_j} = \exp \left( \langle \phi_j, \delta_i \rangle
\right), \qquad j=1, \ldots, n,$$
which we denote by
$$\exp(\Theta_{ji}) := \exp \left( \langle \phi_j, \delta_i \rangle
\right), \qquad j=1, \ldots, n.$$
Here the index $i$ is still fixed. Note that we have no problem to pair
the $\Theta_{ji}$ to the $\Phi_j, j=1, \ldots, n$ since for each $i$
the $A_{ji}$ are paired to the $\alpha_j, j=1, \ldots, n$ through the
Vandermonde-like systems \eqref{VM} and \eqref{VMlike}.

By setting up \eqref{VMlike} for each 
$i=1, \ldots, d-1$ and pairing its solution with $\Phi_j$ in \eqref{Phi},
we obtain for fixed $j=1, \ldots, n$ the linear system of equations
\begin{equation}
\begin{pmatrix} \Delta_1 & \cdots & \Delta_d \\
\delta_{11} & \cdots & \delta_{1d} \\ \vdots & & \vdots \\
\delta_{d-1,1} & \cdots & \delta_{d-1,d} \end{pmatrix} 
\begin{pmatrix} \phi_{j1} \\ \vdots \\ \phi_{jd} \end{pmatrix} =
\begin{pmatrix} \Phi_j \\ \Theta_{j1} \\ \vdots \\ \Theta_{j,d-1} \end{pmatrix}.
\label{IDshift}
\end{equation}
Since the vectors $\Delta$ and $\delta_i, i=1, \ldots, d-1$ are linearly independent,
the coefficient matrix of \eqref{IDshift} is regular and so the individual
$\phi_{ji}, j=1, \ldots, n, i=1, \ldots, d$ can be computed, at the expense
of $2n$ evaluations $F_s$ in \eqref{1stbatch} and $(d-1)n$ evaluations 
$F_{\ell i}$ in \eqref{IDsamples}.

Before we continue we point out that (as is clear from the semantics of
the formulas) we can also denote $\Delta$ as $\delta_0$, 
$F_s$ as $F_{s0}$ and $\Phi_j$ as $\Phi_{j0}$. 

\section{Disentangling collisions}

We now turn our attention to the situation in which the first batch of
samples $F_{s}$ at multiples of the vector $\Delta$ does not reveal all
individual terms because some values $\exp(\Phi_j), j=1 \ldots, n$ collide
and the exponential sum shrinks to $\nu<n$ terms.
For ease of notation, but without loss of generality, we take the 
colliding terms to be successive, for instance: $\exp(\Phi_1) = \ldots =
\exp(\Phi_{h_1}), \exp(\Phi_{h_1+1}) = \ldots = \exp(\Phi_{h_2}), \ldots,
\exp(\Phi_{h_{\nu-1}+1}) = \ldots = \exp(\Phi_n)$.
Assume that with $0 \le s \le 2\nu-1, \nu \le n$
the exponential samples break down into
\begin{multline}
F_s = \sum_{j=1}^\nu \left( \alpha_{h_{j-1}+1} + \ldots + \alpha_{h_j} \right) 
\exp \left(
\langle \phi_{h_j}, s\Delta \rangle \right), \\ \qquad h_0=0, \quad h_j < h_{j+1},
\quad h_\nu = n \label{collision}
\end{multline}
because
$$\exp(\Phi_{h_{j-1}+1}) = \ldots = \exp(\Phi_{h_j}), \qquad j=1, \ldots,
\nu.$$
Since $|\Im\phi_{ji}| < \pi/|\Delta_i|, j=1, \ldots, n, i=1, \ldots, d$,
we actually have
$$\Phi_{h_{j-1}+1} = \cdots = \Phi_{h_j}, \qquad j=1, \ldots, \nu.$$
The Vandermonde system \eqref{VM} now becomes
\begin{multline}
\hskip -0.6truecm \begin{pmatrix} 
\exp(k\Phi_{h_1}) & \cdots & \exp(k\Phi_{h_\nu}) \\ 
\vdots & & \vdots \\
\exp((k+\nu-1)\Phi_{h_1}) & \cdots & \exp((k+\nu-1)\Phi_{h_\nu}) 
\end{pmatrix} \begin{pmatrix} \alpha_1 + \ldots + \alpha_{h_1}
\\ \vdots \\ \alpha_{h_{\nu-1}+1} + \ldots + \alpha_{h_\nu}
\end{pmatrix} \\ = 
\begin{pmatrix} F_{k} \\ \vdots \\
F_{k+\nu-1} \end{pmatrix}, \qquad 0 \le k \le \nu.
\label{VMnu}
\end{multline}
Note that at the same time, the degree of the polynomial $B(z)$ in
\eqref{pol} is only
$\nu$. How this is detected and how the true $n$ is revealed is discussed in
the next section. To proceed we denote
\begin{equation}
A_j := \alpha_{h_{j-1}+1} + \ldots + \alpha_{h_j}, \qquad j=1, \ldots, \nu. \label{Aj}
\end{equation}

To disentangle the collisions in the exponential sum, we need additional
evaluations besides the minimal number $(d+1)n$. At the end of 
Section 5 we also explain how these additional evaluations allow to deal
with the situation where some $A_j=0$.

We start with $i=1$ and the identification shift vector $\delta_1$. 
First we point out how the Vandermonde-like system \eqref{VMlike} 
of Section 3 looks like in
case of such collisions: in the coefficient matrix the value $n$ is
replaced by $\nu$ and in $\Phi_j$ the index $j$ is replaced by $h_j$.
With the collisions in \eqref{collision}, the
unknowns $A_{j1}, j=1, \ldots, \nu$ take the form
$$A_{j1} = \alpha_{h_{j-1}+1} \exp\left( \langle \phi_{h_{j-1}+1}, \delta_1
\rangle \right) +
\ldots + \alpha_{h_j} \exp \left( \langle \phi_{h_j}, \delta_1 \rangle
\right), 
\qquad j=1, \ldots, \nu.$$
In the sequel we denote from here on the additional evaluations $F_{\ell
1}$ mentioned in Section 3 by
$$F_{1\ell 1}:= F_{\ell 1} = f(\kappa_{\ell 1}\Delta + \delta_1), \qquad 
\ell = 1, \ldots, \nu$$
and we add, still with $i=1$, the samples 
$$F_{s \ell 1}:= f(\kappa_{\ell 1}\Delta+s\delta_1), \quad 
s=2, 3, \ldots, 2\max_{1 \le j \le \nu}(h_j-h_{j-1}), \quad 
\ell=1, \ldots, \nu, \quad i=1.$$
The triple index expresses the shift vector multiple in the index $s$, 
the collision into $\nu$ piles of the $\Phi_j$ in the index $\ell$, 
and the identification level in $i$ (which is $i=1$ here).

Since the values of $h_j$ are actually unknown, the addition of samples is
done further and interlaced with singularity checks of some Hankel
matrices, as we explain now. The checks are performed for each collision
or pile $h_j$ and later repeated
for each $i$. Collisions in the space spanned by $\Delta$ may not be
fully disentangled in the space spanned by $\Delta$ and $\delta_1$, but they
are gradually being disentangled as we add independent vectors $\delta_i$ 
until we span the whole space.
At the last stage, when dealing with the full basis
$\Delta, \delta_1, \ldots, \delta_{d-1}$, the true $n$ is revealed 
because in the end all collisions are taken apart, given enough additional samples. 
For the moment we continue
with $i=1$. 

For each $s$ separately, we set up in analogy with \eqref{VMlike}, the 
Vandermonde-like system 
\begin{equation}
\begin{pmatrix} \exp(\kappa_{h_1 1} \Phi_{h_1}) & \exp(\kappa_{h_1 1}\Phi_{h_2}) & 
\cdots & \exp(\kappa_{h_1 1}\Phi_{h_\nu}) \\ \\
\exp(\kappa_{h_2 1}\Phi_{h_1}) & \cdots & & \exp(\kappa_{h_2 1}\Phi_{h_\nu}) \\ 
\vdots & & & \vdots \\ \\ \exp(\kappa_{h_\nu 1} \Phi_{h_1}) & \cdots
& & \exp(\kappa_{h_\nu 1}\Phi_{h_\nu}) \end{pmatrix}
\begin{pmatrix} A_{s11} \\ \vdots \\ A_{s\nu 1} \end{pmatrix} = \begin{pmatrix}
F_{s 1 1} \\ \vdots \\ F_{s \nu 1} \end{pmatrix} \label{disVMlike}
\end{equation}
where
\begin{multline}
A_{sj1} = \alpha_{h_{j-1}+1} \exp\left( \langle \phi_{h_{j-1}+1}, s\delta_1
\rangle \right) +
\ldots + \alpha_{h_j} \exp \left( \langle \phi_{h_j}, s\delta_1 \rangle \right), \\
\qquad j=1, \ldots, \nu. \label{Adis}
\end{multline}
Note that the coefficient matrix is independent of $s$. Also, the former
unknowns $A_j$ and $A_{j1}$ can as well be indexed as $A_{0j1}$
and $A_{1j1}$ respectively, and so \eqref{disVMlike}
and \eqref{Adis} remain valid for $s=0, 1$, which is important for the sequel. 
The values $A_j$ from \eqref{Aj} and $A_{sj1}, s \ge 1$ from \eqref{Adis} are 
actually equidistant samples of the function
\begin{align}
A_{j1}(x)=  & A_{j1}(x_1, \ldots, x_d) \notag \\
=   & \alpha_{h_{j-1}+1} \exp\left( \langle 
\phi_{h_{j-1}+1}, x \rangle \right) +
\ldots + \alpha_{h_j} \exp \left( \langle \phi_{h_j}, x \rangle \right), \notag \\
& \hspace{7.5cm} j=1, \ldots, \nu,
\label{Ax}
\end{align}
taken at $x=s\delta_1, s \ge 0$. 
For each fixed $j=1, \ldots, \nu$ we now put together the Hankel matrix
\begin{equation}
\begin{pmatrix}
A_{0j1} & A_{1j1} & A_{2j1} & A_{3j1} & \ldots \\
A_{1j1} & A_{2j1} & A_{3j1} & \ldots &  \phantom{\vdots} & \\
A_{2j1} & A_{3j1} & A_{4j1} & \ldots & \phantom{\vdots} \\
A_{3j1} & \vdots & \vdots & \\
\vdots &
\end{pmatrix}. \label{disHan}
\end{equation}
Note that in order to enlarge \eqref{disHan} with one row and column for a
particular $j$, one needs to
solve \eqref{disVMlike} for two additional values of $s$, thereby obtaining the
additional $A_{sj1}$ for all $1 \le j \le \nu$.

It is known that the rank of any $(h_j-h_{j-1}+t) \times (h_j-h_{j-1}+t)$ 
submatrix for finite $t \ge 0$ is bounded by $h_j-h_{j-1}$
\cite{Ka.Le:ear:03,Cu.Le:spa:16} since $h_j-h_{j-1}$ equals the number of terms 
in each of the evaluations $A_j, 
A_{sj1}, s \ge 1$. The actual rank $r_j$ of the $(h_j-h_{j-1}) \times
(h_j-h_{j-1})$ submatrix with $A_j$ in the top left corner tells us (with high probability \cite{Ka.Le:ear:03}) how
many of the $h_j-h_{j-1}$ terms in $A_j(x)$ can indeed be separated at
the current level ($i=1$)
where identification shift $\delta_1$ is brought into the
picture.
The value of $r_j$ is 
discovered as one adds samples $F_{s\ell 1}$, 
solves \eqref{disVMlike} and enlarges \eqref{disHan} step by step. 
This explains why we need to add samples $F_{s\ell 1}$ 
until $s$ reaches $2\max_j (h_j-h_{j-1})$ or until for all $j$ the rank
$r_j$ is known. How do we proceed to extract the
coefficients and exponential parameters from \eqref{Ax} and disentangle the
collisions? 

For $j$ fixed, $r_j$ of the individual terms 
$$\alpha_{h_{j-1}+k} \exp(\langle \phi_{h_{j-1}+k}, \delta_1 \rangle ), 
\qquad k=1, \ldots, h_j-h_{j-1}, 1 \le r_j \le h_j-h_{j-1}$$ 
of $A_{j1}(x)$ can be deduced from the samples $A_j, 
A_{sj1}, s\ge 1$ of $A_{j1}(x)$
using one of the Prony-like methods
\cite{Be.Ti:det:88,Sc:mul:86,Ro.Ka:esp:89,Hu.Sa:mat:90} which were already 
mentioned to solve for \eqref{Phi} from \eqref{HA} or \eqref{EV} and 
compute the coefficients from \eqref{VM}.
We remark that when 
$r_j< h_j- h_{j-1}$ then some collisions in $A_{j1}(x)$ still remain indistinguishable in the
space spanned by $\Delta$ and $\delta_1$.

For the sake of completeness we explicitly give the generalized eigenvalue 
problems that
lead to the identification of $1 \le r_j \le h_j-h_{j-1}$ terms in $A_j(x)$:
\begin{multline*}
\begin{pmatrix} A_{1j1} & A_{2j1} & \cdots & A_{r_j, j1} \\
A_{2j1} & A_{3j1} & \cdots & A_{r_j+1, j1} \\
\vdots & & \vdots \\
A_{r_j, j1} & A_{r_j+1, j1} & \cdots & A_{2r_j-1, j1} \\
\end{pmatrix} v  \\ = \lambda
\begin{pmatrix} 
A_{0j1} & A_{1j1} & \cdots & A_{r_j-1, j1} \\
A_{1j1} & A_{2j1} & \cdots & A_{r_j, j1} \\
\vdots & & \vdots \\
A_{r_j-1, j1} & A_{r_j, j1} & \cdots & A_{2r_j-2, j1} \\
\end{pmatrix} v, 
\qquad v \in \cz^{r_j}.
\end{multline*}

After disentangling at $i=1$, at least partially, some of
the collisions, we can update the number of terms in the
exponential model from $\nu$ to $\mu \ge \nu$ and reduce the collisions to
\begin{multline*}
F_s = \sum_{j=1}^{\mu} \left( \alpha_{g_{j-1}+1} + \ldots +
\alpha_{g_j} \right) \exp \left(
\langle \phi_{g_j}, s\Delta \rangle \right), \\ g_0=0, \quad
g_j < g_{j+1}, \quad g_{\mu} = n.
\end{multline*}
It is clear that 
the previous indices $h_j, j=1, \ldots, \nu$ are among the $g_k,
k=1, \ldots, \mu$ but remember that we don't know the values $h_j$ or $g_k$ 
explicitly. We only know that
for some $j$ a collision from index $h_{j-1}+1$ to $h_j,
1 \le j \le \nu$ may have split into separate piles indexed by some
$g_k$ and $g_{k+1}, 1 \le k \le \mu$.
At this moment in the procedure, namely at the completion of
step $i=1$, we have computed
$$\exp(\Theta_{g_k 1}), \qquad \Theta_{g_k 1} := \langle \phi_{g_k}, \delta_1
\rangle, \qquad k= 1, \ldots, \mu.$$
Because $|\Im \langle \phi_j, \delta_1/||\delta_1|| \rangle | <
\pi/||\delta_1||, j=1, \ldots, n$
we in fact obtained all the values
$$\Theta_{g_{k-1}+1, 1} = \cdots = \Theta_{g_k, 1}, \qquad k=1, \ldots, \mu,
\qquad g_0=0, \quad, g_\mu=n$$
which we need later on in combination with the
$$\Phi_{h_{j-1}+1} = \cdots = \Phi_{h_j}, \qquad j=1, \ldots, \nu \le \mu$$
to identify the individual $\phi_{ji}$ as in
\eqref{IDshift}. 

We now explain how to move from $i$ to $i+1$. 
The first thing is to find
proper locations for the samples involving the next identification shift
$\delta_2$.
Some care needs to be taken with respect to the
regularity of the Vandermonde matrices involved.
For $i=2$ we collect
\begin{multline}
F_{s \ell 2}:= f(\kappa_{\ell 2} (\Delta+\delta_1) + s \delta_2), 
\\ s=1, \ldots, 2\max_{1 \le j \le \mu} (g_j-g_{j-1}), \quad 
\ell=1, \ldots, \mu, \quad i=2. \label{collect}
\end{multline}
Let us denote
$$\Omega_{g_j 1} := \Phi_{g_j} + \Theta_{g_j 1}, 
\qquad j=1, \ldots, \mu.$$
Note that the sum is a direct consequence of the choice $\Delta+\delta_1$ in
\eqref{collect}, which is briefly discussed below.
Similarly to \eqref{disVMlike} we write down, for each $s$ separately,
\begin{equation}
\begin{pmatrix} \exp(\kappa_{g_1 2} \Omega_{g_1 1}) & 
\exp(\kappa_{g_1 2} \Omega_{g_2 1}) & 
\cdots & \exp(\kappa_{g_1 2} \Omega_{g_{\mu} 1}) \\ \\
\exp(\kappa_{g_2 2} \Omega_{g_1 1}) & \cdots & & 
\exp(\kappa_{g_2 2} \Omega_{g_{\mu} 1}) \\ 
\vdots & & & \vdots \\ \\ \exp(\kappa_{g_{\mu} 2} \Omega_{g_1 1}) & \cdots
& & \exp(\kappa_{g_{\mu} 2} \Omega_{g_{\mu} 1}) \end{pmatrix}
\begin{pmatrix} A_{s12} \\ \vdots \\ A_{s\mu 2} \end{pmatrix} 
= \begin{pmatrix} F_{s 1 2} \\ \vdots \\ F_{s \mu 2} \end{pmatrix}
\label{i2}
\end{equation}
where
\begin{multline}
A_{sj2} = \alpha_{g_{j-1}+1} \exp\left( \langle \phi_{g_{j-1}+1}, 
s\delta_2 \rangle \right) +
\ldots + \alpha_{g_j} \exp \left( \langle \phi_{g_j}, 
s\delta_2 \rangle \right), \\ j=1, \ldots, \mu. \label{A2sj}
\end{multline}
From here it is clear how to finalize the $i=2$ phase and how to proceed to
the next value of $i$. We point out that instead of the linear combination
$\Delta+\delta_1$ in \eqref{collect}, any linear combination
$c\Delta+e\delta_1$ with $ce\not=0$ that guarantees the regularity of the
coefficient matrix in \eqref{i2} can be used (then the definition of $\Omega_{g_j}$ also
needs to be adapted). This option may be useful as
it allows to control the location of the sample points for numeric purposes
or so.

To round up this section, we
summarize the algorithm that recovers the vectors $\phi_j$ and
coefficients $\alpha_j$ for $j=1, \ldots, n$ in case of possible
collisions of inner products with the chosen directional vectors $\Delta,
\delta_i, i=1, \ldots, d-1$. Before we proceed, we further adapt our
notation. Let
\begin{align*}
\delta_0 &:= \Delta, \\
\nu_{-1} &:= 0, \\
\nu_0 &:= \nu, \\
\nu_1 &:= \mu
\end{align*}
Our first aim is to identify all the inner products 
$\Phi_{ji} = \langle \phi_j, \delta_i \rangle, j=1, \ldots, n, i=0,
\ldots, d-1$, including possible collisions. 
This is done by making use of successively collected
samples, namely
\begin{multline*}
F_{s \ell i} = f\left(\kappa_{\ell i}(\delta_0+ \ldots + \delta_{i-1}) + s
\delta_i \right), \\ s=0, 1, 2, \ldots \quad \ell= 1, 2, \ldots,
\nu_{i-1}, \quad i=0, \ldots, d-1,
\end{multline*} 
where we assume that empty sums equal zero and values in an empty range
need not be specified. The samples are collected by fixing the indices
from the right to the left: at identification level $i$, collision or pile
$\ell$ is being sparsely 
interpolated using the samples collected at shift multiples
$s$. Here $\nu_i$ indicates
the number of non-coinciding inner products at identification level $i$.
Remember that $s$ is running up to twice the number of terms in expression
$A_{\ell i}(x)$ at level $i$ (for $i=1$ this is given in \eqref{Ax} and
it is straightforward 
to imagine how it looks like for general $i$). We remind
the reader that only the evaluation at multiples of $\delta_i, i \ge 0$
needs to follow an equidistant scheme. The values $\kappa_{\ell i}$ need
not be like that. We also mentioned earlier that the sum $\delta_0+ \ldots
+\delta_{i-1}$ can be replaced by another linear combination. The only
crucial element is that the $\delta_i, i \ge 0$ are linearly independent.
The latter will precisely allow us to identify the vector components
$\phi_{ji}, j=1, \ldots, n, i=1, \ldots, d$ from the inner products 
$\langle \phi_j, \delta_i \rangle, j=1, \ldots, n, i=0, \ldots, d-1$
as in \eqref{IDshift}.

\section{Detecting the sparsity}

The minimal number of $(d+1)n$ samples only delivers the parameters
$\alpha_j, \phi_{ji}, j=1, \ldots, n, i=1, \ldots,d$ if the
value of $n$ is somehow known a priori and no collision of values
$\exp(\Phi_j), j=1, \ldots, n$ occurs. In the previous section we described
how to deal with eventual collisions. Here we detail how to detect the value
of $n$ should it not be given. In addition, we analyze how many samples
are needed in the worst case when neither $n$ is known nor the
projections are collision free.

While collecting the samples $F_s=f(s\Delta)$ and building the Hankel
matrices in \eqref{HA} or \eqref{EV}, 
the rank of the Hankel matrix reveals (with high probability \cite{Ka.Le:ear:03}) the
number $\nu$ of terms that do not collide when evaluating in the space spanned by
the vector $\Delta$. To this end we need at least $2\nu+1$ values so that we
can compose the $(\nu+1)\times (\nu+1)$ Hankel matrix 
$$\begin{pmatrix}
F_0 & \ldots & F_\nu \\
\vdots & & \vdots \\
F_\nu & \ldots & F_{2\nu}
\end{pmatrix}$$
and conclude that it
is singular \cite{He:app:74,Ba.Gr:pad:81,Ka.Le:ear:03}.

From $\nu$ and \eqref{HA} or \eqref{EV}
we proceed to collect the samples $F_{1\ell 1}$ ($s=1$) and
$F_{2\ell 1}$ ($s=2$), another $2\nu$ in total ($\ell=1, 2, \ldots, \nu$). 
If all $2 \times 2$ Hankel 
matrices of the form \eqref{disHan} are singular, then every collision
remains indistinguishable (unless the zero determinant was an unfortunate
coincidence \cite{Ka.Le:ear:03})
also in the space spanned by $\Delta$ and $\delta_1$. 
However, if for some $j$ the 
$2 \times 2$ matrix \eqref{disHan} is regular, then we have to
proceed to the next values for $s$ ($s=3, 4$), collect another $2\nu$ values in
total, and find out how many terms actually can be revealed in the space spanned by
$\Delta$ and $\delta_1$. We proceed until we find no
larger matrices of the form \eqref{disHan} that are regular. Only after
working ourselves through all regular matrices of the form \eqref{disHan} with
$\delta_1$ ($i=1$) we can update $\nu$ to $\mu \ge \nu$. 

And then we bring the next identification shift vector $\delta_2$
in the picture. We collect the samples $F_{s\ell 2}$ ($s=1, 2$)
as in \eqref{collect}
and compose matrices similar to \eqref{disHan} but now with the last index
in the $A_{sj1}$ replaced by $i=2$ and with $A_{sj2}$ defined as in 
\eqref{A2sj}.
The inspection of the Hankel matrices
containing the values computed for $A_{sj2}$ is identical to the procedure
described in the previous paragraph for $i=1$. If required, as before, 
we add more samples for larger values of $s$.

Finally, by the time we reach $i=d-1$ we can update the number
of terms to the true value for $n$. Now how many samples has this cost us?
When $n$ is known a priori and
we do not run into collisions or cancellations, which with high
probability do not occur, the algorithm presented in Section 3 uses only
$$(d+1)n$$
samples. Next, we look at the situation where collisions occur and Section
4 is put to work (how to deal with possible
cancellations is dealt with at the end of this section). Also the sparsity
$n$ is not given.
The $A_j$ and $\Phi_{h_j}$ with $j=1, \ldots, \nu$ in \eqref{collision}
are retrieved from $O(\nu)$ samples where $\nu \le n$. In $A_j, 1 \le j
\le \nu$ there are
$h_j-h_{j-1}$ terms colliding, where each $h_j-h_{j-1} \le n-\nu+1$. To
disentangle the terms in $A_j$ we need $O(h_j-h_{j-1})$ samples and so we
need at most $O(\nu (n-\nu+1))$ samples to disentangle all $A_j, j=1,
\ldots, \nu$. Note that we have overestimated each $h_j-h_{j-1}$ by
$n-\nu+1$, while if one $h_j-h_{j-1}=n-\nu+1$, all others equal 1. The
procedure is repeated when working with the identifications shifts
$\delta_1, \ldots, \delta_{d-1}$, leading us to a grand total of
\begin{equation}
O\left( (d+1) \max_{1 \le \nu \le n} \nu(n-\nu+1) \right).
\label{datacomplexity}
\end{equation}

Remains to discuss the issue of a vanishing $A_{sji}$. For simplicity, but
without loss of generality, we discuss the situation where one of the
coefficients $A_j$ given by \eqref{Aj} vanishes, in other words
$A_j=A_{j1}(0)=0$ with $A_{j1}(x)$ given by \eqref{Ax}. So besides
encountering a collision, the result of the collision is now also zero.

If some $A_j=0$ then the rank of the matrices in \eqref{EV} is less than
$\nu$ and will not reveal the correct value for $\nu$. 
Of course, the accidental cancellation of a coefficient $A_j$
happens only with very small probability.
It suffices either to probe $f(x_1,
\ldots, x_d)$ along another (random) choice for the vector $\Delta$
\cite{Zippel1979,Ka.Le:ear:03}, or if one 
absolutely wants to extract the information $\langle \phi_j, \Delta
\rangle$ for the originally chosen $\Delta$, 
to probe $f(x_1, \ldots, x_d)$ along one or more (random)
parallel shifts of $\Delta$, as in 
\begin{equation}
F_s:=f(s\Delta_1+k\epsilon, \ldots, s\Delta_d+k\epsilon), \qquad
s=0, \ldots, n, \qquad k=1,2, \ldots \label{cancel}
\end{equation}
Such a shift affects the coefficient $A_j$ in that it changes from
$A_j(0)$ to 
$$A_j(k\epsilon) = 
\alpha_{h_{j-1}+1} \exp\left( \langle \phi_{h_{j-1}+1}, k\epsilon \rangle 
\right) + \ldots + \alpha_{h_j} \exp \left( \langle \phi_{h_j}, k\epsilon 
\rangle \right).$$
The rank of the matrices in \eqref{EV} when filled with the values in
\eqref{cancel} either confirms the already computed rank $\nu$ or reveals
a higher and more probably correct rank $\nu$.
The random probing or parallel translation can be added to every step
$i=0, \ldots, d-1$ in the procedure when selecting $\delta_0=\Delta,
\delta_1, \ldots, \delta_{d-1}$ without impacting our data usage analysis
in \eqref{datacomplexity}.

All the above is now illustrated with an example 
in which we take the reader through the entire process, 
first collision-free, then including collision disentanglement.

\section{Numerical illustration}

We take $d=2$, write $u:=x_1, v:=x_2, x=(u,v)^t$ 
and consider the exponential sum
$$f(u,v) = \sum_{j=1}^4 \alpha_j \exp(\langle \phi_j, x \rangle)$$
with
\begin{equation*}
\begin{aligned}
\phi_1 &= (-0.5, 1 + \im 2 \pi \times 0.5), \\
\phi_2 &= (0.1 + \im 2 \pi \times 3.4, 1.5 + \im 2 \pi \times 5.2), \\
\phi_3 &= (0.1 + \im 2 \pi \times 3.4, -0.5 + \im 2 \pi \times 12.6), \\
\phi_4 &= (-2.5 + \im 2 \pi \times 23.2, -10 + \im 2 \pi \times 82.3),
\end{aligned} \qquad
\begin{aligned}
\alpha_1 &= 1.7 \exp(\im 2 \pi/10), \\
\alpha_2 &= 1.1 \exp(\im 2 \pi/20), \\
\alpha_3 &= 0.9, \\
\alpha_4 &= 9.2 \exp(\im 2 \pi/2).
\end{aligned}
\end{equation*}
When outputting numerical results for this small scale example, we round
all values to 4 significant digits (all relative errors are less than $5
\times 10^{-4}$). The numerical effect of the choice of
the vectors $\Delta$ and $\delta_i$ throughout the process, and that of the 
underlying one-dimensional Prony-like
method in use, is beyond the scope of this paper and will be the subject
of further investigations.

First we show the simple case described in the Sections 2 and 3, where the
number of terms $n=4$ is known up front and no collisions of the inner
products in the samples occur. Of course, the latter is hard to predict in
practice.

We take $\Delta=(0.01, 0.01)$ and $\delta_1=(-0.01, 0.01)$. Using 8
equidistant evaluations at $x=s\Delta, s=0, \ldots, 7$, we obtain 
from \eqref{EV} the
values of $\exp(\Phi_j)$ and can deduce the $\Phi_j, j=1,
\ldots, 4$ because $|\Im \phi_{ji}| < \pi/|\Delta_i|$:
\begin{align*}
\Phi_1 &= <\phi_1, \Delta> \approx  0.005000 +
0.03142\im, \\
\Phi_2 &= <\phi_2, \Delta> \approx  0.01600 +
0.5404\im, \\
\Phi_3 &= <\phi_3, \Delta> \approx  -0.004000 +
1.005\im, \\
\Phi_4 &= <\phi_4, \Delta> \approx -0.1250 +
0.3456\im.
\end{align*}
We obtain the coefficients $\alpha_j, j=1, \ldots, 4$
from \eqref{VM}:
\begin{align*}
\alpha_1 &\approx 
1.700 \exp(\im 2 \pi \times 0.1000), \\
\alpha_2 &\approx 
1.100 \exp(\im 2 \pi \times 0.05000), \\
\alpha_3 &\approx 
0.9000 \\ 
\alpha_4 &\approx 
9.200 \exp(\im 2 \pi \times 0.5000).
\end{align*}
From 4 additional evaluations along
the identification shift $\delta_1$, we obtain the values of
$\exp(\Theta_{11})$, $\exp(\Theta_{21})$, $\exp(\Theta_{31})$,
$\exp(\Theta_{41})$ from \eqref{VMlike}. Their exponents are the projections of the vectors $\phi_j$
along $\delta_1$:
\begin{align*}
\Theta_{11} &= <\phi_1, \delta_1> \approx 0.01500 +
0.03142\im, \\
\Theta_{21} &= <\phi_2, \delta_1> \approx 0.01400 +
0.1131\im, \\
\Theta_{31} &= <\phi_3, \delta_1> \approx -0.006000 +
0.5781\im, \\
\Theta_{41} &= <\phi_4, \delta_1> \approx -0.07500 +
3.713\im.
\end{align*}
We finally obtain the values of $\phi_j = (\phi_{j1},\phi_{j2})^t$ 
by solving for each $j=1, \ldots, 4$
$$\begin{pmatrix} \Delta_1 & \Delta_2 \\ \delta_{11} & \delta_{12}
\end{pmatrix} \begin{pmatrix} \phi_{j1} \\ \phi_{j2} \end{pmatrix} =
\begin{pmatrix} \Phi_j \\ \Theta_{j1} \end{pmatrix}.$$
This leads to the following numerical approximations for the $\phi_j$:
\begin{align*}
\phi_1
& \approx (-0.5000, 1.000 + \im 2 \pi \times 0.5000), \\
\phi_2 
& \approx (0.1000+ \im 2 \pi \times 3.400, 1.500 + \im 2 \pi \times 5.200),\\
\phi_3 
& \approx (0.1000 + \im 2 \pi \times 3.400, -0.5000 + \im 2 \pi \times 12.60), \\
\phi_4 
& \approx (-2.500 + \im 2 \pi \times 23.20, -10.00 + \im 2 \pi \times 82.30).
\end{align*}
So far we have used 12 samples in total, which indeed equals $(d+1)n$.
Next we deal with the situation in which neither $n$ is known, nor the
assumption of the non-collision holds.

One additional evaluation in the first batch, at $x=8\Delta$,
would ideally (meaning that the numerical rank
is easy to detect) and with high probability (meaning that we don't
accidentally hit a root of the determinant)
have revealed that $n=4$, still under the assumption 
that no collisions occur at the inner products. But let us instead move to
other directions $\Delta$ and $\delta_1$ that get us in trouble because of
colliding inner products. 

Take $\Delta=(0.03, 0)$ and $\delta_1=(0, 0.01)$.
The projections of $\phi_2$ and $\phi_3$ along $\Delta$ clearly coincide.
After 7 evaluations at $x=s\Delta, s=0, \ldots, 6$ we found that $\nu_0=3$
and we obtain from \eqref{VMnu}
that (without actually knowing the values of the $h_j$ which we list
only to help the reader follow the example):
\begin{align*}
\Phi_{h_1} &= <\phi_1, \Delta> 
\approx -0.01500, \\ 
\Phi_{h_2} &= <\phi_3, \Delta> 
\approx 0.003000 + 0.6409\im, \\
\Phi_{h_3} &= <\phi_4, \Delta> 
\approx  -0.07500 + 4.373\im.
\end{align*}
We proceed without knowing $n$ and without knowing whether and where some
collisions have occurred. But we know, since $d=2$, that after adding an
independent shift vector $\delta_1$, all terms will have revealed themselves.

So we add evaluations $F_{s \ell 1} = f(\kappa_{\ell 1}\Delta+s\delta_1)$
with $\ell=1, 2, 3$ and $s=1, 2, \ldots$ For simplicity we choose
$\kappa_{\ell 1} = \ell-1$.
With $\ell=1$ and $s=1, 2$ we find that the matrix
$$\begin{pmatrix} 
A_1 & A_{111} \\ A_{111} & A_{211}
\end{pmatrix},$$
where the $A_{sji}$ are computed from \eqref{Adis},
has rank 1 and so $h_1=1=g_1$.
With $\ell=2$ and $s=1, 2, 3, 4$ we find that the matrix
$$\begin{pmatrix} 
A_2 & A_{121} & A_{221} \\ A_{121} & A_{221} & A_{321} \\
A_{221} & A_{321} & A_{421}
\end{pmatrix}$$
has rank 2. This indicates with high probability that there are 2 terms coinciding at
$\Phi_{h_2}$ (hence $h_2=3$ and $g_2=2, g_3=3$).
Remember that in order to obtain $A_{s21}, 1 \le s \le 4$, we need to solve
\eqref{disVMlike} which involves the samples $F_{sj1}, 1\le j \le 3$.
Hence, continuing the sampling for $\ell=2$ drags along $\ell=1, 3$ at the
same time. In other words, we are now spending $3 \times 4$ samples
for $\ell=1,2,3$ rather than only 4 samples for $\ell=2$.

We now reveal $\langle \phi_2, \delta_1 \rangle$ and $\langle \phi_3,
\delta_1 \rangle$ by solving the generalized eigenvalue
problem
$$\begin{pmatrix} A_{121} & A_{221} \\ A_{221} & A_{321}
\end{pmatrix} v = \lambda \begin{pmatrix} 
A_2 & A_{121} \\ A_{121} & A_{221} \end{pmatrix} v.$$
With $\ell=3$ and $s=1,2$ we find the same conclusion as with $\ell=1$,
now for
$$\begin{pmatrix} 
A_3 & A_{131} \\ A_{131} & A_{231}
\end{pmatrix},$$
and so $\nu_1=4$ with $h_3=4, g_4=4$. 

At the expense of a total of $(2 \times 3 + 1)+ 3 \times 4 = 19$ 
evaluations, we find that $n=4$ and we can
identify all $\phi_{ji}$ and $\alpha_j$ for $j=1, \ldots, 4$ and $i=1,2$.

\section{Conclusion}

In 1795 the French scientist G.~de Prony showed that a univariate 
linear combination of $n$ exponential terms with unknown real but mutually distinct
exponents
could be fitted uniquely to $2n$ data samples. His result solves the $d=1$ case of this
paper. 
The current paper is the first of its kind where this result is proven to
hold for
general $d>1$: a multivariate linear combination of $n$ exponential terms with 
unknown inner product exponents can, under mild conditions,
be fitted using only $(d+1)n$ data.

\section*{Acknowledgements}

This work was partially supported by a Research Grant of the FWO-Flanders (Flemish Science Foundation). 

\bibliographystyle{spmpsci}      

          
\end{document}